\documentclass{article}
\usepackage{latexsym,amssymb,graphicx,tikz,adjustbox,tikz-cd,amsmath}
\input amssym.def \input amssym

\def\dj{d\kern-0.4em\char"16\kern-0.1em}
\def\Dj{\mbox{\raise0.3ex\hbox{-}\kern-0.4em D}}
\newtheorem{te}{Theorem}[section]
\newtheorem{de}[te]{Definition}
\newtheorem{lm}[te]{Lemma}

\newtheorem{ex}[te]{Example}

\def\bN{\mathbb{N}}
\def\dokaz{\noindent{\bf Proof. }}
\def\kraj{\hfill $\Box$ \par \vspace*{2mm} }
\usepackage[draft,nomargin,inline]{fixme}
\fxsetface{inline}{\itshape}
\fxsetface{env}{\itshape}
\fxusetheme{color}

\newcommand{\newinf}{\mathop{\mathrm{inf}\vphantom{\mathrm{sup}}}}

\begin{document}
\begin{center}
           {\huge \bf Measures of string similarities based on the Hamming distance}
\end{center}
\begin{center}
{\small \bf Bojan Nikoli\' c, Boris  \v Sobot}\\[2mm]
{\small  Faculty of Natural Sciences and Mathematics, University of Banja Luka,\\
Mladena Stojanovi\'ca 2, 78000 Banja Luka, Bosnia and Herzegovina,\\
e-mail: bojan.nikolic@pmf.unibl.org\\

Faculty of Sciences, University of Novi Sad,\\
Trg Dositeja Obradovi\'ca 4, 21000 Novi Sad, Serbia\\
e-mail: sobot@dmi.uns.ac.rs}
\end{center}
\begin{abstract} \noindent
In this paper we consider measures of similarity between two sets of strings built up using the Hamming distance and tools of persistence homology as a basis.
First we describe the construction of the \v Cech filtration adjoined to the set of strings, the persistence module corresponding to this filtration and its barcode structure. Using these means, we
introduce a novel similarity measure for two sets of strings, based on a comparison of bars within their barcodes of the same dimension. Our idea is to look for
a comparison that will take under consideration not only the overlap of bars, but also ensure that observed bars are qualitatively matched, in the sense that they represent
similar homological features. To make this idea happen, we developed a method called the separation of simplex radii technique.
\vspace{1mm}\\

{\sl 2020 Mathematics Subject Classification}: 05E45, 55N31,  62R40 \\

{\sl Key words and phrases}: Hamming distance, string similarity measure, \v Cech filtration, persistence module, barcode, bottleneck distance,
separation of simplex radii technique
\end{abstract}

\section{Introduction}


A string is a finite sequence over  a (usually finite) alphabet. We will consider strings on $n$-element alphabets and assume (without loss of generality) that strings are over the alphabet $\mathbb{N}_n=\{1,2,\dots,n\}$. By $S(n,l)$ we denote the set of strings of length $l$ over this alphabet. We also consider a string $s\in S(n,l)$ as a function $s:\mathbb{N}_l\rightarrow\mathbb{N}_n$ and denote its $i$-th character by $s(i)$.\\

The Hamming distance between two strings of equal length is the number of positions at which the corresponding symbols are different. More precisely,
the Hamming distance on $S(n,l)$ is defined as follows: for $s=a_1a_2\dots a_l$ and $t=b_1b_2\dots b_l$:
$$d_H(s,t):=|\{i\in\{1,2,\dots,l\}:a_i\neq b_i\}|.$$
This type of edit distance was introduced by R. W. Hamming in his seminal paper \cite{H}, and has applications in several disciplines, including information theory, coding theory,
cryptography, and bioinformatics. \\

Let $A$ and $B$ be subsets of $S(n,l)$ with same cardinality $m\geqslant 1$. In the case $m=1$, the Hamming distance
between the element of set $A$ and the element of set $B$ can be used as a measure of dissimilarity of this sets. In the case $m>1$, a measure of dissimilarity of sets
$A$ and $B$ in the metric space $(S(n,l),d_H)$ can be defined via the Hausdorff distance between these sets:
$$ D_H(A,B):=\max\{\sup_{a\in A}\newinf_{b\in B}d_H(a,b),\sup_{b\in B}\newinf_{a\in A} d_H(a,b)\}.     $$

The Hausdorff distance $D_H(A,B)$ is "one-dimensional" in its nature and does not consider the internal structures of sets $A$ and $B$. This disadvantage motivates us to
consider various types of connectivity (or lack thereof) that exist between the elements of sets $A$ and $B$ separately. In this way, it is possible to define a
similarity measure between these sets which would compare their connectivity classes within the same dimension. In this article, we will accomplish this by using tools from simplicial homology as well as its multiscale version known as persistent homology (\cite{EH},\cite{ZC}). \\

In recent years, simplicial homology and persistent homology have played a central role in Topological data analysis (TDA), a branch of applied mathematics which analyzes topological information from high-dimensional datasets. Simplicial homology studies the sequence of homology groups of a simplicial complex. Loosely speaking, objects of a homology group
are holes of a specific dimension that reside in the observed simplicial complex. The most notable type of simplicial complex is the \v Cech complex, which is defined as the nerve
of the cover of balls of a fixed radius around each point of a given set of points.
In our setting, for  arbitrary $r\geqslant 0$ and $A\subseteq S(n,l)$, the
\v Cech complex $\mathcal C_A^{(r)}$ is a simplicial complex consisting of all nonempty $\sigma\subseteq A$ such that the closed balls of radius $r$ with centers in $\sigma$ have a non-empty common intersection.

For a given integer $k\geqslant 0$ and \v Cech complex $\mathcal C_A^{(r)}$, the homology group of dimension $k$ will be denoted by $H_k\left(\mathcal C_A^{(r)}\right)$.
Elements of this group are $k-$dimensional homological classes, i.e. $k-$cycles on $\mathcal C_A^{(r)}$ which are not boundaries. The group $H_k\left(\mathcal C_A^{(r)}\right)$ captures $k-$dimensional topological features when the \v Cech complex $\mathcal C_A^{(r)}$ is observed  with resolution $r$. In  most cases, we don't have enough information
which would enable us to choose the "optimal" resolution $r$. Hence, it is useful to observe the \v Cech filtration, the family of \v Cech complexes $\{\mathcal C_A^{(r)}:r\geq 0\}$ obtained by varying resolution (level) $r$ in the definition of $\mathcal C_A^{(r)}$. Clearly, for $r_1<r_2$ holds $\mathcal C_A^{(r_1)}\subseteq \mathcal C_A^{(r_2)}$, and, since $A$ is a finite set of strings, "blowing up" resolution would lead to a level $r_t$ of filtration such that simplicial complex $\mathcal C_A^{(r_t)}$ is the full complex, that is, it contains every nonempty subset of $A$. Therefore, all \v Cech filtrations that we consider in this paper have a finite number of levels, i.e., they have a form $\{\mathcal C_A^{(r_1)},\mathcal C_A^{(r_2)},\dots, \mathcal C_A^{(r_t)}\}$,
for some $0\leqslant r_1<r_2<\dots<r_t$. We will call filtrations $\{\mathcal C_A^{(r_i)}:i\in\{1,2,\dots,t\}\}$ and $\{\mathcal C_B^{(r_i)}:i\in\{1,2,\dots,t\}\}$ isomorphic if there is a bijection $f:A\rightarrow B$ such that, for each $i\in\{1,2,\dots,t\}$, $\sigma\in\mathcal C_A^{(r_i)}$ if and only if $f[\sigma]\in\mathcal C_B^{(r_i)}$.
An automorphism of the metric space $(S(n,l),d_H)$ mapping $A$ to $B$ is called  a $d_H(A\rightarrow B)-$isomorphism.\\

Persistent homology keeps track of the evolution of homological classes throughout the levels of a given filtration. More precisely, for each dimension $k$, the persistence module
\begin{equation} \label{pers.module}
H_k\left(\mathcal C_A^{(r_1)}\right)\hookrightarrow H_k\left(\mathcal C_A^{(r_2)}\right)\hookrightarrow\dots\hookrightarrow H_k\left(\mathcal C_A^{(r_t)}\right)
\end{equation}
contains information on the complete lifespan of every $k-$dimensional homology class, from the level of filtration when they are first formed (born) to the level of filtration when they become boundaries, and hence trivial (die). In this way, we obtain the persistence interval $[\mbox{birth}(\gamma),\mbox{death}(\gamma))$, for every homological class $\gamma$. In \cite{ZC},
it was shown that persistence module (\ref{pers.module}) has a unique decomposition into a multiset of persistence intervals. This multiset is called the barcode of dimension
$k$ and is denoted by $BC_k$. Visually, a barcode $BC_k$ can be represented by a multiset of horizontal lines whose endpoints correspond to the birth-death pairs of $k-$homological
classes. We will abuse notation and use $BC_k(A)$ to denote the $k-$dimensional barcode of the persistence module corresponding to the \v Cech filtration adjoined to the subset
$A\subseteq S(n,l)$. \\

Comparing persistence barcodes is appealing due to their simple numerical nature. In the context of persistent homology, the most notable measure of
comparison is the bottleneck distance. The idea behind this distance is to observe all possible bijections (matchings) between two multisets of barcode lines,
such that every line of "significant length" from one barcode is paired with a unique line of similar length and endpoints from the other barcode, and vice versa. The
bottleneck distance between two barcodes is then defined as an infimum of the set of significant lengths for which described matching can be done  (see the next section
for the precise definition). The most important property of the bottleneck distance is its stability, in the sense that "small" changes in the structure of
the persistence module lead to small amount of changes in the corresponding barcode \cite{SEH}. One notable matching between two barcodes which enables proof
of this stability is induced matching introduced in \cite{BL}. \\

The results we present in this paper are focused on studying the similarity of two subsets $A,B\subseteq S(n,l)$ of the same cardinality $m$. The aforementioned similarity will be expressed through the appropriate matching, which would generate a measure of similarity between subsets $A$ and $B$.\\

Under this framework, the main contributions of this paper are as follows:

\begin{itemize}
\item We introduce notions of generalized strings and generalized Hamming distance. These concepts allow us to develop a novel simplices radii separation technique.
This technique is based on constructing a bijection which maps a subset $A\subseteq S(n,l)$ to an appropriate subset $A'$ of generalized strings so that the following two useful properties are satisfied:

1. All barcode lines of $BC_k(A')$, for $k\geqslant 1$, have unique birth-death endpoints. As a consequence, we can propose a fairly simple matching between barcodes $BC_k(A')$
 and $BC_k(B')$.

2. Changes in the structure of the persistence module (\ref{pers.module}) that occur after applying this bijection are strictly controlled.

\item  The simplices radii separation technique enables us to consider a new sort of barcode matching based on the idea of cycle registration. This matching
allows us to induce a novel similarity measure betwen two subsets $A,B\subseteq S(n,l)$.
\end{itemize}
The rest of this paper is organized as follows. Section \ref{sec:preliminaries} sets up basic notions and properties of simplicial homology and persistent homology.
Section \ref{sec:isom.filtr.} introduces \v Cech filtration adjoined to the set $A\subseteq S(n,l)$. Also in this section,
notions of generalized strings and generalized Hamming distance are introduced. Section \ref{sec:barcode} gives details about barcodes for a subset of strings. Furthermore, a
bijection between  a subset $A\subseteq S(n,l)$ and an appropriate subset $A'$ of generalized strings is provided. This bijection yields a barcode $BC_k(A')$, which is
"close enough" to the barcode $BC_k(A)$ and has useful "nonaligned" setup of its lines. This property will be used in order to define appropriate "hybrid" matching between
barcodes of two subsets of strings. The last section presents conclusions and plans for future work.\\

\section{Preliminaries} \label{sec:preliminaries}

In this section, we briefly recall the basic concepts of simplicial homology and persistence homology. For a more in-depth examination, see for example \cite{ZC},
\cite{SEH}, \cite{BL}, \cite{EH2} and \cite{PRSZ}.\\

A \emph{simplicial complex} $\mathcal K$ is a pair $(K,\Sigma)$, where $K$ is a nonempty set and $\Sigma$ is a finite collection of nonempty subsets of $K$ called \emph{simplices}, such that $\tau\subseteq\sigma\in\Sigma$ implies $\tau\in\Sigma$. The simplex $\sigma$ with elements $v_0,\dots,v_k$ is denoted by $[v_0,\dots,v_k]$ instead of
$\{v_0,\dots,v_k\}$. The \emph{dimension of a simplex} $\sigma$ is $dim\,\sigma = |\sigma|-1$ and the \emph{dimension of the complex} is the maximum dimension of all of its simplices. \emph{Full complex} is a simplicial complex $(K,P(K)\setminus\{\emptyset\})$. If $\tau\subseteq\sigma$, then $\tau$ is a \emph{face} of $\sigma$.
The \emph{vertex set} of the complex $\mathcal K$ is the collection of all elements $v\in K$ such that $v\in\sigma$, for some simplex $\sigma\in\Sigma$, and is denoted by $Vert(\mathcal K)$.
A \emph{subcomplex} $\mathcal L$ of the complex $\mathcal K=(K,\Sigma)$ is a simplicial complex whose simplices form a subfamily of $\Sigma$. For simplicial complexes
$\mathcal K=(K,\Sigma_K)$ and $\mathcal L=(L,\Sigma_L)$, a mapping $f : Vert(\mathcal K)\rightarrow Vert(\mathcal L)$ such that $\sigma\in \Sigma_K$ if and only if $f[\sigma]\in\Sigma_L$ is called a \emph{simplicial mapping}. Two simplicial complexes are \emph{isomorphic} if there is a simplicial bijection between these complexes. \\

Let $(X,d)$ be a metric space. For $r>0$ and $x\in X$, let $B(x,r)=\{y\in X:d(x,y)\leq r\}$ be the closed ball of radius $r$ around $x$. If $K\subseteq X$ is a finite set,
for every $r>0$ the \emph{\v Cech complex} $\mathcal C_K^{(r)}$ is the simplicial complex $(K,\{A\in P(K)\setminus\{\emptyset\}:\bigcap_{x\in A}B(x,r)\ne\emptyset\})$. For $r_1<r_2$, $\mathcal C_K^{(r_1)}$ is a subcomplex of $\mathcal C_K^{(r_2)}$, which we write (informally) as
$\mathcal C_K^{(r_1)}\subseteq \mathcal C_K^{(r_2)}$. \\

Let $\mathcal K$ be a simplicial complex and $k$ a dimension. A \emph{$k$-chain} is a formal sum $c=\sum a_i\sigma_i$, where the $\sigma_i$ are the $k$-simplices
(i.e., simplices with the dimension $k$) in $\mathcal K$  and the $a_i$ are coefficients from the field $\mathbb Z_2$. Addition of two $k-$chains is defined
componentwise, i.e. if $c_1=\sum a_i\sigma_i$ and $c_2=\sum b_i\sigma_i$, then $c_1+c_2=\sum(a_i+b_i)\sigma_i$. For every dimension $k$, the $k$-chains together with the
addition operation form the \emph{group of $k$-chains} denoted as $C_k(\mathcal K)$. The \emph{boundary} of the $k$-simplex $\sigma=[v_0,v_1,\dots,v_k]$ is the sum of its $(k-1)$-dimensional faces, i.e.,
$\partial_k\sigma=\sum_{j=0}^k [v_0,v_1,\dots,\hat{v_j},\dots,v_k],$
where the hat indicates that $v_j$ is omitted. For an arbitrary $k$-chain, its boundary is the sum of the boundaries of its simplices.
A \emph{$k$-cycle} $c$ is a $k$-chain with empty boundary, $\partial_k c = 0$. Since $\partial$ commutes with addition, we have a group of $k$-cycles, denoted as
$Z_k(\mathcal K)= ker \partial_k$.
A \emph{$k$-boundary} $c$ is a $k$-chain that is the boundary of a $(k + 1)$-chain, $c = \partial_{k+1} d$, with $d\in C_{k+1}$. Since $\partial$ commutes with addition,
we have a group of $k-$boundaries, denoted by $B_k(\mathcal K)= im \partial_{k+1}$.
The \emph{$k$-th homology group} is the $k$-th cycle group modulo the $k$-th boundary group, $H_k(\mathcal K) = Z_k(\mathcal K)/B_k(\mathcal K)$.
Each element of $H_k = H_k(\mathcal K)$ is obtained by adding all $k$-boundaries to a given
$k$-cycle, $c + B_k$, with $c\in Z_k$, and this class is referred as a \emph{homology class}. Nontrivial homology classes (for $c\neq 0$) depict cycles that are not boundaries of any chain of simplices of appropriate dimension. In the language of a geometric realization of the given complex, these cycles represent "holes" of suitable dimensions.
Every simplicial map $f$ between simplicial complexes $\mathcal K$ and $\mathcal L$
can be extended to the \emph{induced homomorphism on homology} $f_k:H_k(\mathcal K)\rightarrow H_k(\mathcal L)$, which maps cycles to cycles and boundaries to boundaries.
The most notable type of induced homomorphism $f_k:H_k(\mathcal K)\rightarrow H_k(\mathcal L)$ occurs in the case when $\mathcal K$ is a subcomplex of $\mathcal L$, i.e.\ when the simplicial map $f$ is an inclusion.  \\

A \emph{filtration} of the simplicial complex $\mathcal K$ is a collection $\{\mathcal K^{(i)}:i\in\{0,1,\dots,t\}\}$ of increasing subcomplexes of $\mathcal K$:
\begin{equation} \label{filtracija}
\emptyset =\mathcal K^{(0)}\subseteq\mathcal K^{(1)}\subseteq\dots\subseteq\mathcal K^{(t)}=\mathcal K.
\end{equation}
If the complex $\mathcal K$ contains $u$ simplices, then a filtration of this complex can be understood as a construction of $\mathcal K$ by adding $t\leqslant u$ chunks of
its simplices, one at a time. For the purpose of an enumeration of complexes in the given filtration, the set $\{0,1,\dotsc,t\}$ can be replaced with an arbitrary
set $\{r_0,r_1,\dotsc,r_t\}$, such that $r_0<r_1<\dots<r_t$.

\begin{ex} For a given \v Cech complex $\mathcal C_K^{(r)}$, every collection $\{r_0,r_1,r_2,\dots,r_t\}$, such that $r_0<0\leqslant r_1<r_2<\dots<r_t=r$, determines a filtration
$$\emptyset=\mathcal C_K^{(r_0)}\subseteq\mathcal C_K^{(r_1)}\subseteq\dots\subseteq\mathcal C_K^{(r_t)}=\mathcal C_K^{(r)},  $$
The \v Cech complex $\mathcal C_K^{(r_i)}$ can be interpreted as the "state" of the complex $\mathcal C_K^{(r)}$ at a resolution level
$r_i\leqslant r$. Thus, radius $r_i$ is also referred to as the level of the filtration. It is worth noting that the values $r_i$ can be chosen in such a way that each stage of the construction has exactly one representative, more precisely: that for every $r'>0$ there is a unique $r_i\leq r'$ such that $\mathcal C_K^{(r_i)}=\mathcal C_K^{(r')}$. In such case we call this the \v Cech filtration.
\end{ex}

For every $i\leqslant j$ and each dimension $k$, we have the induced homomorphism $f_k^{i,j} : H_k(\mathcal K^{(r_i)})\rightarrow H_k(\mathcal K^{(r_j)})$ generated by
the inclusion map $\mathcal K^{(r_i)}\hookrightarrow\mathcal K^{(r_j)}$. The filtration thus corresponds to a sequence of homology groups connected by homomorphisms:
\begin{equation} \label{istr.modul}
0 =H_k(\mathcal K^{(r_0)})\stackrel{f_k^{0,1}}{\rightarrow} H_k(\mathcal K^{(r_1)})\stackrel{f_k^{1,2}}{\rightarrow}\dots\stackrel{f_k^{t-1,t}}{\rightarrow} H_k(\mathcal K^{(r_t)})=H_k(\mathcal K),
\end{equation}
again, one for each dimension $k$. The sequence (\ref{istr.modul}) is also called the \emph{persistence module} and is denoted by $PM_k(\mathcal K)$. As we go from $\mathcal K^{(r_{i-1})}$ to $\mathcal K^{(r_i)}$,
we gain new homology classes and we lose some when they become trivial or merge with each other. We collect the classes that are born at or before a given threshold
and die after another threshold in groups. The \emph{$k$-th persistent homology groups} are the images of the homomorphisms induced by inclusion, $H_k^{i,j} = im f_k^{i,j}$,
for $0\leqslant i\leqslant j\leqslant t$.
Note that $H_k^{i,i} = H_k(\mathcal K^{(r_i)})$. The persistent homology groups consist of the homology classes of $\mathcal K^{(r_i)}$ that are still alive at $\mathcal K^{(r_j)}$
or, more formally, $H_k^{i,j} = Z_k(\mathcal K^{(r_i)})/\left(B_k(\mathcal K^{(r_j)})\cap Z_k(\mathcal K^{(r_i)})\right)$. We have such a group for each dimension $k$ and each
index pair $i\leqslant j$.
A homology class $\gamma\in H_k(\mathcal K^{(r_i)})$ is \emph{born at} $\mathcal K^{(r_i)}$ if $\gamma\notin H_k^{i-1,i}$. Furthermore, if $\gamma$ is born at $\mathcal K^{(r_i)}$, then it \emph{dies entering} $\mathcal K^{(r_j)}$, if it merges with an older class as we go from $\mathcal K^{(r_{j-1})}$ to $\mathcal K^{(r_j)}$, that is, $f_k^{i,j-1}(\gamma)\notin H_k^{i-1,j-1}$, but $f_k^{i,j}(\gamma)\in H_k^{i-1,j}$.  A \emph{positive simplex} is a simplex with property that its addition in some level of filtration leads to the birth of a new homology class. Similarly, a \emph{negative simplex} is a simplex with property that its addition in some level of filtration leads to the death of an existing homology class. If $\gamma$ is born at $\mathcal K^{(r_i)}$ and dies entering $\mathcal K^{(r_j)}$, then the interval $[i,j)$ is called the \emph{persistence interval} of the homology class $\gamma$. The length of this interval is called the \emph{persistence of the homology class} $\gamma$ and is denoted by $pers(\gamma)$. If a homology class $\gamma$ is born at $\mathcal K^{(r_i)}$ but never dies, then the interval $[i,+\infty)$ is the persistence interval of this class and we set $pers(\gamma)=\infty$. Persistence intervals keep track of the lifespan of all homology classes in the process of passing through the observed persistence module. A filtration having the property that, at every level of the filtration, the homology changes allowed are either the creation of a single new cycle or the termination of a single existing cycle, is called a \emph{Morse filtration}. Essentially, all persistence intervals of the persistence module corresponding to a Morse filtration have different endpoints.
The notion of persistence module can also be defined for a sequence of vector spaces that are not necessarily homology groups. \\

A \emph{morphism} between persistence modules $PM_k(\mathcal K)$ and $PM_k(\mathcal L)$ given by
\begin{align*}
&0 =H_k(\mathcal K^{(r_0)})\stackrel{f_k^{0,1}}{\rightarrow} H_k(\mathcal K^{(r_1)})\stackrel{f_k^{1,2}}{\rightarrow}\dots\stackrel{f_k^{t-1,t}}{\rightarrow} H_k(\mathcal K^{(r_t)})=H_k(\mathcal K),\\
&0 =H_k(\mathcal L^{(r_0)})\stackrel{g_k^{0,1}}{\rightarrow} H_k(\mathcal L^{(r_1)})\stackrel{g_k^{1,2}}{\rightarrow}\dots\stackrel{g_k^{t-1,t}}{\rightarrow} H_k(\mathcal L^{(r_t)})
=H_k(\mathcal L),
\end{align*}
is a collection $h=\left\{h_i:H_k(\mathcal K^{(r_i)})\rightarrow H_k(\mathcal L^{(r_i)}):i\in\{0,1\dots,t\}\right\}$ of homomorphisms, such that, for every $i<j$, the
following diagram is commutative:

\begin{center}
\begin{tikzcd}
    H_k(\mathcal K^{(r_i)}) \arrow[d,"{h_i}"]\arrow[r,"{f_k^{i,j}}"] & H_k(\mathcal K^{(r_j)}) \arrow[d,"{h_j}"] \\
    H_k(\mathcal L^{(r_i)}) \arrow[r,"{g_k^{i,j}}"] & H_k(\mathcal L^{(r_j)})
\end{tikzcd}
\end{center}
A morphism $h$ connecting persistence modules $PM_k(\mathcal K)$ and $PM_k(\mathcal L)$ is also denoted by $h:PM_k(\mathcal K)\Rightarrow PM_k(\mathcal L)$.
Specially, if every map in its collection is a bijection, then $h$ is an \emph{isomorphism} and, in this case, persistence modules $PM_k(\mathcal K)$ and
$PM_k(\mathcal L)$ are \emph{isomorphic persistence modules}. If $h:PM_k(\mathcal K)\Rightarrow PM_k(\mathcal L)$ is a morphism, then a persistence module given by
$$ 0 =h_0\left[H_k(\mathcal K^{(r_0)})\right]\stackrel{g_k^{0,1}}{\rightarrow} h_1\left[H_k(\mathcal K^{(r_1)})\right]\stackrel{g_k^{1,2}}{\rightarrow}\dots\stackrel{g_k^{t-1,t}}{\rightarrow} h_t\left[H_k(\mathcal K^{(r_t)})\right]
=h_t\left[H_k(\mathcal K)\right],  $$
is called the \emph{image of morphism} $h$, and is denoted by $im(h)$. For a $\delta>0$, the \emph{$\delta-$shifted persistence module} $PM_k(\mathcal K)(\delta)$ is obtained by
"shifting" levels of the module $PM_k(\mathcal K)$ to the left by $\delta$, i.e., at the $i-$th level of this module is the homology group $H_k(\mathcal K^{(r_i+\delta)})$
and the induced homomorphism connecting $i-$th and $j-$th level of this module is equal to $f_k^{r_{i'},r_{j'}}$, where $i',j'$ are such that $r_{i'}\leq r_i+\delta<r_{i'+1}$ and $r_{j'}\leq r_j+\delta<r_{j'+1}$. The \emph{$\delta-$shifted morphism} $h^{\delta}$ between
persistence modules $PM_k(\mathcal K)$ and $PM_k(\mathcal K)(\delta)$ is a morphism given by the collection $\{h_i^{\delta}:i\in\{0,1\dots,t\}\}$, such that, for every $i$,
$h_i^{\delta}:=f_k^{i,i+\delta}$. Also, for a morphism $h:PM_k(\mathcal K)\Rightarrow PM_k(\mathcal L)$, the morphism between their corresponding $\delta-$shifted modules
is denoted by $h(\delta)$.
For a $\delta>0$, persistence modules $PM_k(\mathcal K)$ and
$PM_k(\mathcal L)$ are \emph{$\delta-$interleaved} if there exist two morphisms $F:PM_k(\mathcal K)\Rightarrow PM_k(\mathcal L)(\delta)$ and
$G:PM_k(\mathcal L)\Rightarrow PM_k(\mathcal K)(\delta)$, such that, for every $i$, the following diagrams are commutative:
$$\adjustbox{scale=0.75,center}{
\begin{tikzcd}
      H_k(\mathcal K^{(r_i)}) \arrow[dr,"{F_i}"] \arrow[r,"{h_i^{\delta}}"] 
      & H_k(\mathcal K^{(r_i+\delta)})\arrow[r,"{h_{i+\delta}^{\delta}}"] & H_k(\mathcal K^{(r_i+2\delta)}) \\
      & H_k(\mathcal L^{(r_i+\delta)}) \arrow[ur,"{G_{i+\delta}}"] &
\end{tikzcd}
\begin{tikzcd}
      & H_k(\mathcal K^{(r_i+\delta)}) \arrow[dr,"{F_{i+\delta}}"] & \\
      H_k(\mathcal L^{(r_i)}) \arrow[ur,"{G_i}"] \arrow[r,swap,"{h_i^{\delta}}"] 
      & H_k(\mathcal L^{(r_i+\delta)})\arrow[r,swap,"{h_{i+\delta}^{\delta}}"] & H_k(\mathcal L^{(r_i+2\delta)})
      \end{tikzcd}
}
$$
The \emph{interleaving distance} between persistence modules $PM_k(\mathcal K)$ and $PM_k(\mathcal L)$ is defined as the infimum of the set of all $\delta>0$
for which modules $PM_k(\mathcal K)$ and $PM_k(\mathcal L)$ are $\delta-$interleaved. This distance is denoted by $d_{INT}$ and it can be proven that it is an
extended pseudo-metric on the set of all persistence modules. \\

A \emph{barcode} is a finite multiset of intervals, i.e.\ a finite collection of intervals with given multiplicities. The intervals in a barcode are also called
\emph{bars}. One notable example of a barcode is the multiset of all persistence intervals corresponding to the persistence module $PM_k(\mathcal K)$. This barcode is denoted by $BC_k\left(PM_k(\mathcal K)\right)$. For every bar $[b,d)$ in this barcode, we can define the \emph{interval persistence module} $I[b,d)$:\\

\adjustbox{scale=0.75,center}{
\begin{tikzcd}
      0 \arrow[r,"{0}"] & \dots \arrow[r,"0"] &  0\arrow[r,"{0}"] & \mathbb Z_2 \arrow[r,"{id_{\mathbb Z_2}}"] & \mathbb Z_2 \arrow[r,"{id_{\mathbb Z_2}}"] & \dots \arrow[r,"{0}"] & 0 \arrow[r,"{0}"] & 0 \arrow[r,"{0}"] & \dots\\
      0 \arrow[u] & \dots & i<b \arrow[u] & i=b \arrow[u] & i\in(b,d) \arrow[u] & \dots & i=d \arrow[u] & i>d \arrow[u] & \dots
\end{tikzcd}
}\\

We have the representation $ PM_k(\mathcal K)=\bigoplus_{[b_i,d_i)} I[b_i,d_i)^{m_i},  $
where $m_i$ is the multiplicity of the persistence interval $[b_i,d_i)$, which belongs to the persistence module $PM_k(\mathcal K)$. This result is known as the
\emph{Normal Form Theorem for Persistence Modules} and was first proved in \cite{B}. As a consequence, every persistence module is completely determined by the structure of
bars in its barcode. \\

Given an interval $I = [b,d)$, denote by $I^{\delta} = [b-\delta,d+\delta)$ the interval obtained by "stretching" $I$ by $\delta$ on both sides.
Let $BC_k(\cdot)$ be a barcode. For $\varepsilon>0$, denote by $BC_k^{\varepsilon}(\cdot)$ the set of all bars from $BC_k(\cdot)$ with length greater than $\varepsilon$.
A \emph{matching} between two finite multisets $X$ and $Y$ is a relation $\mu\subseteq X\times Y$, such that $\mu:X'\rightarrow Y'$ is a bijection between some $X'\subseteq X$ and $Y'\subseteq Y$. In this case, $coim(\mu)=X'$, $im(\mu)=Y'$, and the elements of $X'$ and $Y'$ are \emph{matched}. If an element appears in the multiset several times, we treat its different copies separately, e.g. it could happen that only some of its copies are matched. If $\mathcal K$ and $\mathcal L$ are filtered complexes and $k\geqslant 0$ a dimension, then a
\emph{$\delta-$matching} between barcodes $BC_k\left(PM_k(\mathcal K)\right)$ and $BC_k\left(PM_k(\mathcal L)\right)$ is a matching $\mu$ which satisfies the following properties:
\begin{align*}
& 1.\,\, BC_k^{2\delta}\left(PM_k(\mathcal K)\right)\subseteq coim(\mu),\\
& 2.\,\, BC_k^{2\delta}\left(PM_k(\mathcal L)\right)\subseteq im(\mu),\\
& 3.\,\, \mbox{If } \mu(I)=J,\mbox{ then } I\subseteq J^{\delta}\mbox{ and } J\subseteq I^{\delta}.
\end{align*}
The \emph{bottleneck distance}, $d_{BOT}\left(BC_k\left(PM_k(\mathcal K)\right), BC_k\left(PM_k(\mathcal L)\right)\right)$ is defined to be the infimum over all $\delta\geqslant 0$
for which there is a $\delta$-matching between barcodes $BC_k\left(PM_k(\mathcal K)\right)$ and $BC_k\left(PM_k(\mathcal L)\right)$. The fundamental property of
the bottleneck distance is stated in the next theorem, the proof of which can be found in \cite{SEH} or \cite{CSGGO}.

\begin{te}(The Isometry Theorem) For persistence modules $PM_k(\mathcal K)$ and $PM_k(\mathcal L)$ holds
$$  d_{INT}\left(PM_k(\mathcal K),PM_k(\mathcal L)\right)=d_{BOT}\left(BC_k(PM_k(\mathcal K)), BC_k(PM_k(\mathcal L)\right). $$
\end{te}

The claim $d_{BOT}\left(BC_k(PM_k(\mathcal K)), BC_k(PM_k(\mathcal L)\right)\leqslant d_{INT}\left(PM_k(\mathcal K),PM_k(\mathcal L)\right)$ is also known as \emph{The Stability Theorem}. Intuitively, this theorem guarantees that "little tweaks" in the structure of the persistence module do not produce significant changes in the structure of bars within
the barcode. \\

\section{Filtration of a set of strings} \label{sec:isom.filtr.}

Let us recall that $S(n,l)$ denotes the set of all strings of length $l$ over the alphabet $\mathbb N_n=\{1,2,\dots,n\}$. For strings $s=a_1a_2\dots a_l$ and
$t=b_1b_2\dots b_l$ in $S(n,l)$, we observe the Hamming distance between them defined by:
$$d_H(s,t):=|\{i\in\{1,2,\dots,l\}:a_i\neq b_i\}|.$$

Also, remember that, for a given $r\geqslant 0$, $B(s,r)=\{t\in S(n,l):d_H(s,t)\leq r\}$ denotes the closed ball of radius $r$ around the element $s$ in the metric space
$(S(n,l),d_H)$. In this section, we will describe the construction of the \v Cech filtration adjoined to the subset $A\subseteq S(n,l)$. Then, we will generalize this procedure in the case of the set of generalized strings.

\subsection{The \v Cech filtration adjoined to $A\subseteq S(n,l)$}

Let $A\subseteq S(n,l)$ be an arbitrary nonempty set of strings. For an arbitrary $r\geqslant 0$, we can consider the \v Cech complex
$\mathcal C_A^{(r)}$, whose simplices are all subsets $\sigma\subseteq A$ with the property $\bigcap_{s\in\sigma} B(s,r)\ne\emptyset$.
Since $A$ is a finite set, there is a minimal terminal radius $r_t\geqslant 0$, such that $\mathcal C_A^{(r)}$ is the full complex for every $r\geqslant r_t$. We are going
to consider filtration of the full complex $\mathcal C_A:=\mathcal C_A^{(r_t)}$ that formalizes the idea of describing all "stepping stones" in the process of building
this complex from the initial complex $\mathcal C_A^{(0)}=\{[s]:s\in A\}$. At the first step, we find the smallest value $r_1>0$ with the property
$\mathcal C_A^{(0)}\subsetneq \mathcal C_A^{(r_1)}$. Then, we find the smallest value $r_2>r_1$ with the property $\mathcal C_A^{(r_1)}\subsetneq \mathcal C_A^{(r_2)}$.
Continuing with this process, we eventually come to the the last step $\mathcal C_A^{(r_{t-1})}\subsetneq \mathcal C_A^{(r_t)}=\mathcal C_A$. In this step, all simplices
which were "missing"  in the complex $\mathcal C_A^{(r_{t-1})}$ are added, finishing the construction of the $\mathcal C_A$.

\begin{de}
The filtration $\mathcal C_A^{(0)}\subsetneq\mathcal C_A^{(r_1)}\subsetneq\dots\subsetneq\mathcal C_A^{(r_t)}$ obtained in the previous construction is called
the filtration adjoined to the subset $A\subseteq S(n,l)$ and $\{r_1,\dots,r_t\}$ is also referred to as the set of levels of this filtration.
\end{de}

We remark that the discrete nature of the Hamming distance implies that all levels $r_1,r_2,\dotsc,r_t$ are positive integers.

\begin{de} For a simplex $\sigma\subseteq A$, the smallest $r_{\sigma}\geqslant 0$ with the property $\bigcap_{s\in\sigma} B(s,r_{\sigma})\ne\emptyset$,
is called the radius of $\sigma$. In that case, an arbitrary element $c\in\bigcap_{s\in\sigma} B(s,r_{\sigma})$ is referred to as a center of $\sigma$.
\end{de}

From the previous definition it follows that any radius $r_{\sigma}$ of a simplex $\sigma\in\mathcal C_A$ necessarily has to be one of the levels of the filtration
adjoined to the set $A$. The converse is also true: for any level $r_i$ of the filtration adjoined to the set $A$, the complex $\mathcal C_A^{(r_i)}$ contains some simplex
$\sigma_0$ which is not in the complex $\mathcal C_A^{(r_{i-1})}$, meaning that the radius of this simplex is equal to $r_i$. Note that, unlike the radius of the simplex,
the center of the simplex need not be unique.

\begin{ex} \label{primjer2.1.}Let $A=\{\underbrace{12244131}_{s_1},\underbrace{22223443}_{s_2},\underbrace{32143431}_{s_3},\underbrace{14443214}_{s_4},\underbrace{22134222}_{s_5}\}
\subset S(4,8)$. In order to obtain the \v Cech filtration adjoined to this set, it is sufficient to find the collection of ordered pairs of the form $(\sigma,r_{\sigma})$, where $\sigma\in\mathcal C_A$ and $r_{\sigma}$ is the radius of the simplex $\sigma$:
\begin{align*}
&(\underbrace{[s_1]}_{\sigma_1},0),(\underbrace{[s_2]}_{\sigma_2},0), (\underbrace{[s_3]}_{\sigma_3},0), (\underbrace{[s_4]}_{\sigma_4},0), (\underbrace{[s_5]}_{\sigma_5},0);\\
&(\underbrace{[s_1,s_3]}_{\sigma_6},2);\\
&(\underbrace{[s_1,s_2]}_{\sigma_7},3), (\underbrace{[s_2,s_3]}_{\sigma_8},3), (\underbrace{[s_1,s_2,s_3]}_{\sigma_9},3), (\underbrace{[s_1,s_4]}_{\sigma_{10}},3), (\underbrace{[s_3,s_4]}_{\sigma_{11}},3), (\underbrace{[s_1,s_3,s_4]}_{\sigma_{12}},3),\\
&(\underbrace{[s_1,s_5]}_{\sigma_{13}},3), (\underbrace{[s_2,s_5]}_{\sigma_{14}},3), (\underbrace{[s_3,s_5]}_{\sigma_{15}},3), (\underbrace{[s_1,s_3,s_5]}_{\sigma_{16}},3);\\
&(\underbrace{[s_2,s_4],4}_{\sigma_{17}}), (\underbrace{[s_1,s_2,s_4]}_{\sigma_{18}},4), (\underbrace{[s_2,s_3,s_4]}_{\sigma_{19}},4), (\underbrace{[s_1,s_2,s_3,s_4]}_{\sigma_{20}},4), (\underbrace{[s_1,s_2,s_5]}_{\sigma_{21}},4),\\
&(\underbrace{[s_2,s_3,s_5]}_{\sigma_{22}},4), (\underbrace{[s_1,s_2,s_3,s_5]}_{\sigma_{23}},4), (\underbrace{[s_4,s_5]}_{\sigma_{24}},4), (\underbrace{[s_2,s_4,s_5]}_{\sigma_{25}},4), (\underbrace{[s_3,s_4,s_5]}_{\sigma_{26}},4),\\
&(\underbrace{[s_2,s_3,s_4,s_5]}_{\sigma_{27}},4), (\underbrace{[s_1,s_4,s_5]}_{\sigma_{28}},4);\\
&(\underbrace{[s_1,s_2,s_4,s_5]}_{\sigma_{29}},5), (\underbrace{[s_1,s_3,s_4,s_5]}_{\sigma_{30}},5), (\underbrace{[s_1,s_2,s_3,s_4,s_5]}_{\sigma_{31}},5).
\end{align*}

The radius of the simplex $\sigma_6=[12244131,32143431]$ is $2$, but its center is not unique, e.g.\ the strings $32244431$ and $12143131$ are both centers of this simplex.

From the previous characterization, we derive the required \v Cech filtration:
$$\mathcal C_A^{(0)}\subsetneq\mathcal C_A^{(2)}\subsetneq\mathcal C_A^{(3)}\subsetneq\mathcal C_A^{(4)}\subsetneq\mathcal C_A^{(5)}=\mathcal C_A,$$
where
\begin{align*}
&\mathcal C_A^{(0)}=\{\sigma_i:1\leqslant i\leqslant 5\},\\
&\mathcal C_A^{(2)}=\mathcal C_A^{(0)}\cup\{\sigma_6\},\\
&\mathcal C_A^{(3)}=\mathcal C_A^{(2)}\cup\{\sigma_i:7\leqslant i\leqslant 16\},\\
&\mathcal C_A^{(4)}=\mathcal C_A^{(3)}\cup\{\sigma_i:17\leqslant i\leqslant 28\},\\
&\mathcal C_A^{(5)}=\mathcal C_A^{(4)}\cup\{\sigma_i:29\leqslant i\leqslant 31\}=P(A)\setminus\{\emptyset\}.
\end{align*}

\end{ex}

\begin{de}
Let $A,B\subseteq S(n,l)$ be nonempty subsets of strings for which the filtration adjoined to $A$ and the filtration adjoined to $B$ both have the
identical set of levels $\{r_1,r_2,\dots,r_t\}$. These filtrations are called isomorphic if there is a bijection $f:A\rightarrow B$ (referred to as a filtration isomorphism),
such that for every simplex $\sigma\in\mathcal C_A$ and every $i\in\{1,2,\dots,t\}$ holds: $\sigma\in\mathcal C_A^{(r_i)}$ if and only if $f[\sigma]\in\mathcal C_B^{(r_i)}$.
\end{de}

%

In addition to isomorphism of filtrations, we introduce a somewhat stronger notion.

\begin{de} For subsets $A,B\subseteq S(n,l)$ of the same cardinality, an automorphism of $(S(n,l),d_H)$ mapping $A$ to $B$ is called a $d_H(A\rightarrow B)$-isomorphism. Subsets
$A,B\subseteq S(n,l)$ for which there is a $d_H(A\rightarrow B)-$isomorphism are called $d_H-$isomorphic sets.

\end{de}

It is obvious that $d_H-$isomorphic subsets have isomorphic adjoined filtrations. The converse is generally not true, as we shall see in the following example.

\begin{ex}
Let $l=5$, $n=3$ and $m=3$. Take $s_1=11113$, $s_2=22223$, $s_3=33333$, $s_4=33122$. It is easily checked that the sets $S_1=\{s_1,s_2,s_3\}$ and $S_2=\{s_1,s_2,s_4\}$ yield isomorphic filtrations. Namely, beside vertices, the complex $\mathcal C_{S_1}^{(2)}$ contains 1-simplices $[s_1,s_2]$, $[s_1,s_3]$ and $[s_2,s_3]$,
while $\mathcal C_{S_2}^{(2)}$ contains 1-simplices $[s_1,s_2]$, $[s_1,s_4]$ and $[s_2,s_4]$. Both
$\mathcal C_{S_1}^{(3)}$ and $\mathcal C_{S_2}^{(3)}$ are full complexes, so $\mathcal C_{S_1}^{(2)}$ and $\mathcal C_{S_2}^{(2)}$ are the only nontrivial subcomplexes. This means that the mapping $g:S_1\rightarrow S_2$, defined by $f(s_i)=s_i$, $i\in\{1,2\}$, $f(s_3)=s_4$, is a filtration isomorphism.

On the other hand, there is no $d_H(S_1\rightarrow S_2)$-isomorphism. To show that,
notice that every
$d_H(S_1\rightarrow S_2)$-isomorphism preserves the number $|\{x(i):x\in S_1\}|$ of different letters at some position. Since strings $s_1,s_2,s_3$ all end with the letter $3$,
assumption that $S_1$ and $S_2$ are $d_H-$isomorphic sets would lead to the conclusion that there is a position such that all strings $s_1,s_2,s_4$ have the same letter on that position. However, it is easy to check that this is not the case.
\end{ex}

\subsection{Generalized strings} \label{sub:gen.strings}

The preceding example demonstrates how easy it is to create sets $S_1,S_2\subseteq S(n,l)$, for $n>2$, which are not $d_H$-isomorphic but have isomorphic adjoined filtrations.  In order to reduce the number of such examples, we need to generalize the notion of a string.

\begin{de}
A generalized string of length $l$ over the alphabet $\bN_n$ is a function $s:\bN_l\rightarrow F_n$, where $F_n$ is the set of functions $f:\bN_n\rightarrow[0,1]$ 
such that $\sum_{i=1}^nf(i)=1$. We will denote the set of such generalized strings by $S'(n,l)$, and the image of $i\in\bN_l$ by $s$ will be denoted by $s[i]$. The generalized Hamming distance between $s,t\in S'(n,l)$ is defined by
$$d_{GH}(s,t)=\sum_{i=1}^l\left(1-\sum_{j=1}^n\min\left\{s[i](j),t[i](j)\right\}\right).$$
\end{de}

The distance $d_{GH}(s,t)$ measures the overlapping in functions $s[i]$ and $t[i]$. Every string $s=a_1a_2\dots a_l\in S(n,l)$ can be identified with a generalized string $s$, where $s[i]$ is the function mapping $a_i$ to $1$, and all other letters to $0$, for all $i\in\bN_l$. Using this convention, it is easy to check that $d_{GH}(s,t)=d_H(s,t)$ holds for
arbitrary strings $s,t\in S(n,l)$, so the restriction of $d_{GH}$ to $S(n,l)$ is the "usual" Hamming distance $d_H$.\\

All concepts that we considered in the case of a set $A\subseteq S(n,l)$ (the full complex $\mathcal C_A$, the filtration adjoined to the set $A$, the barcode $BC_k(A)$, etc.) can be introduced analogously in the case of a finite set $A'\subseteq S'(n,l)$. Of course, the diferrence is that we now use the distance $d_{GH}$ instead of $d_H$.

\begin{de}
Let $\mathcal C_{A}$ be the full complex for a finite set $A\subseteq S'(n,l)$. For $\sigma\in \mathcal C_{A}$, the value
$$r(\sigma)=\min\{r:(\exists x\in S'(n,l))(\forall y\in\sigma)d_{GH}(x,y)\leq r\}$$
is called the radius of $\sigma$. A generalized string $c$ such that $d_{GH}(c,y)\leq r(\sigma)$ for all $y\in\sigma$ is called a center of $\sigma$.
\end{de}

\begin{lm}
For every $\sigma\in \mathcal C_{A}$, the minimum in the definition of $r(\sigma)$ exists.
\end{lm}

\dokaz As the set $[0,1]$ with the usual topology is compact, the product space $[0,1]^{\bN_n}$ is also compact. The subspace $S'(n,l)=\{f\in[0,1]^{\bN_n}:\sum_{i=1}^nf(i)=1\}$ is closed, so it is compact as well. If we define a function $\psi_\sigma:S'(n,l)\rightarrow\mathbb{R}$ with
$$\psi_\sigma(x)=\max\{d_{GH}(x,y):y\in\sigma\},$$
it is clearly continuous, so it reaches its minimum on $S'(n,l)$, and that is exactly $r(\sigma)$.\kraj

\begin{ex}

Let $s_1=111112$, $s_2=111113$, $s_3=222221$ and $s_4=333331$. For the set $\sigma=\{s_1,s_2,s_3,s_4\}\subseteq S(3,6)$ one center is the generalized string $b\in S'(3,6)$ given by
\begin{eqnarray*}
b[i]:\left(\begin{array}{lll}
1 & 2 & 3\\
\frac 7{15} & \frac 4{15} & \frac 4{15}
\end{array}\right) & (\mbox{for }1\leq i\leq 5) &
b[6]:\left(\begin{array}{lll}
1 & 2 & 3\\
1 & 0 & 0
\end{array}\right).
\end{eqnarray*}
 Namely, for all $1\leq k\leq 4$ $d_{GH}(b,s_k)=5\cdot\frac 8{15}+1=5\cdot\frac{11}{15}=\frac {11}3$. Note that, for any $c\in S'(3,6)$ and any $1\leq i\leq 5$, $\sum_{k=2}^4(1-\sum_{j=1}^3\min\{c[i](j),s_k[i](j)\})=2$. For $i=6$ we have $\sum_{k=2}^4(1-\sum_{j=1}^3\min\{c[6](j),s_k[6](j)\})\geq 1$, obtaining the minimum only for $c[6]=b[6]$. Thus $\sum_{k=2}^4d_{GH}(c,s_k)\geq 11$, and $r(\sigma)\geq\frac{11}3$.\\

More centers can be obtained by moving weights between first five positions, for example,
\begin{eqnarray*}
b'[1]:\left(\begin{array}{lll}
1 & 2 & 3\\
\frac 5{15} & \frac 5{15} & \frac 5{15}
\end{array}\right) & &
b'[2]:\left(\begin{array}{lll}
1 & 2 & 3\\
\frac 9{15} & \frac 3{15} & \frac 3{15}
\end{array}\right)\\
b'[i]:\left(\begin{array}{lll}
1 & 2 & 3\\
\frac 7{15} & \frac 4{15} & \frac 4{15}
\end{array}\right) & (\mbox{for }3\leq i\leq 5) &
b'[6]:\left(\begin{array}{lll}
1 & 2 & 3\\
1 & 0 & 0
\end{array}\right).
\end{eqnarray*}
\end{ex}

\section{A similarity measure based on comparison of barcodes} \label{sec:barcode}

In this section, we define the barcode associated with a given set of strings. We will use this barcode as an indicator of homological features that appear in the "universe" of the
filtration adjoined to the observed set of strings. Loosely speaking, bars in a barcode represent the evolution of "holes" of appropriate dimension. Every $k-$dimensional hole,
for $k\geqslant 1$, expresses high dimensional "connectivity issue" that exist for some subfamily of strings in some parts of the filtration. Thus, barcodes can be exploited to
measure the discrepancy between connectivity classes of two sets of strings.
 The main goal of this section is to introduce a novel similarity measure for two sets of strings, which would be based on comparison of bars within their
barcodes of the same dimension. This comparison will take under consideration not only the overlap of bars but also ensure that observed bars are "qualitatively" matched, in the
sense that they represent similar homological features.

\subsection{Barcode associated to a set $A\subseteq S(n,l)$}

Let $A\subseteq S(n,l)$ be an arbitrary nonempty set of strings and $\mathcal C_A^{(0)}\subsetneq\mathcal C_A^{(r_1)}\subsetneq\dots\subsetneq\mathcal C_A^{(r_t)}$ the
filtration adjoined to this set. For a fixed dimension $k\geqslant 0$, this filtration generates the persistence module $PM_k(\mathcal C_A)$ given by:
$$
H_k(\mathcal C_A^{(0)})\stackrel{f_k^{0,r_1}}{\rightarrow} H_k(\mathcal C_A^{(r_1)})\stackrel{f_k^{r_1,r_2}}{\rightarrow}\dots\stackrel{f_k^{r_{t-1},r_t}}{\rightarrow} H_k(\mathcal C_A^{(r_t)})=H_k(\mathcal C_A),
$$
where homomorphisms $f_k^{r_i,r_{i+1}}$ are induced by inclusions $\mathcal C_A^{(r_i)}\hookrightarrow \mathcal C_A^{(r_{i+1})}$.

\begin{de} For an arbitrary $k\geqslant 0$, the barcode of the persistence module $PM_k(\mathcal C_A)$ is also referred to as the $k-$dimensional barcode associated
with set $A\subseteq S(n,l)$ and will be denoted by $BC_k(A)$.
\end{de}

Since $A\subseteq S(n,l)$ is a finite set, the barcode $BC_k(A)$ contains no barcode lines (bars) for any $k>|A|-1$.
Barcode $BC_0(A)$ has exactly $|A|$ bars. Each of them shows the evolution of a connected component while moving through the filtration. Since the full complex $\mathcal C_A$ contains only one connection component, we conclude that barcode $BC_0(A)$ has only one infinite bar. For $k\geqslant 1$, every $k-$dimensional hole must eventually be closed
at some level of the filtration, meaning that all bars belonging to the barcode $BC_k(A)$ must have finite lengths.

\begin{ex}
Let us examine barcodes associated with the set $A\subset S(4,8)$ from Example \ref{primjer2.1.}. There are $|A|=5$ bars in the barcode $BC_0(A)$. Since
$[s_1,s_3]\in\mathcal C_A^{(2)}$, two components from $BC_0(A)$ are merged at this level, leaving $4$ bars to persist until the next level of the filtration. Excluding simplices $[s_2,s_4]$ and $[s_4,s_5]$, all other $1-$simplices belong to the complex $\mathcal C_A^{(3)}$, implying that, after this level, there is only one connected component and, consequently, only one bar (with infinite persistence). The barcode $BC_1(A)$ contains only one bar. This bar depicts the persistence of the only nontrivial $1-$cycle $[s_1,s_2]+[s_2,s_5]+[s_1,s_5]$. This cycle is born in the complex $\mathcal C_A^{(3)}$ and dies in the next complex of the filtration,  since $r\left([s_1,s_2,s_5]\right)=4$. Similarly, the barcode $BC_2(A)$ contains only one bar
corresponding to $2-$cycle $[s_1,s_2,s_4]+[s_1,s_2,s_5]+[s_1,s_4,s_5]+[s_2,s_4,s_5]$, which is created in the complex $\mathcal C_A^{(4)}$ and closed down in the full complex $\mathcal C_A^{(5)}$ (see Figure (\ref{fig:barkod dijagram})).

\begin{figure}[h]
\includegraphics[width=1\textwidth]{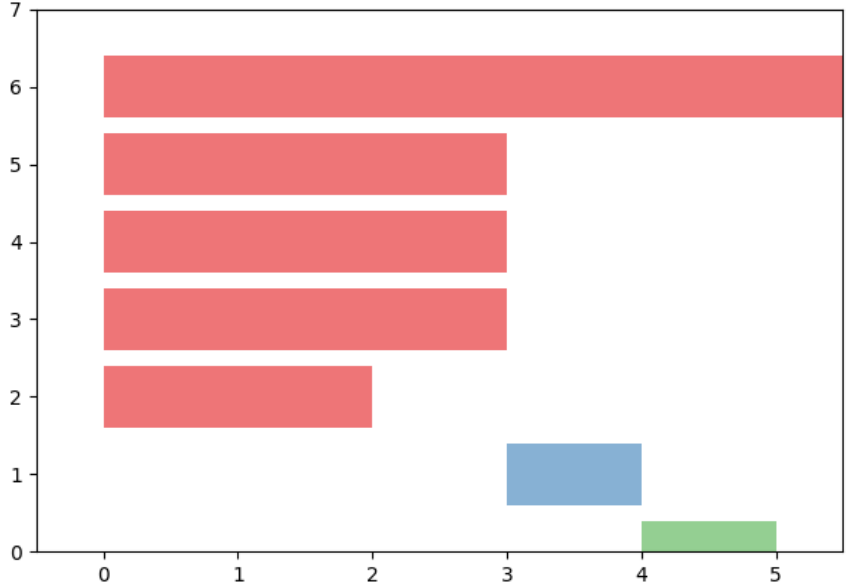}
\centering
\caption{Barcode associated to the given set $A$. Red lines depict bars from the barcode $BC_0(A)$. The blue line depicts the only bar from the barcode $BC_1(A)$ and
the green line depicts the only bar from the barcode $BC_2(A)$.  }
\label{fig:barkod dijagram}
\end{figure}
\end{ex}

Now it is a good time to elaborate on how we are going to use barcodes $BC_k(A)$ and $BC_k(B)$ in order to introduce a measure of (dis)similarity of these sets. Let us suppose
that sets $A,B\subseteq S(n,l)$ have the same cardinality $m\geqslant 2$.

Both of the barcodes $BC_0(A)$ nad $BC_0(B)$ contain $m$ bars, where $m-1$  of them are finite-length bars (all having
$0$ as a left endpoint) and one bar is an infinite-length bar. Two infinite-length bars from these barcodes are perfectly matched and thus can be ignored. The remaining lines in
the $BC_0(A)$ and $BC_0(B)$ can be enumerated in the form  $\left\{[0,l_i^A):i\in\{1,2,\dots,m-1\}\right\}$, $0<l_1^A\leqslant l_2^A\leqslant\dots\leqslant l_{m-1}^A$, and $\left\{[0,l_i^B):i\in\{1,2,\dots,m-1\}\right\}$, $0<l_1^B\leqslant l_2^B\leqslant\dots\leqslant l_{m-1}^B$, respectively. We can match bars $[0,l_i^A)$ and $[0,l_i^B)$ and look for
maximal difference between $l_i^A$ and $l_i^B$.

In the case of dimension $k\geqslant 1$, as a measure of dissimilarity we could use the bottleneck distance  $d_{BOT}\left(BC_k(A),BC_k(B)\right)$. However,
instead of matching the bars exclusively by means of their "best suited" overlaps, we will also try to investigate the possibility of matching bars at a qualitatively higher level.
We will conduct this examination by using the idea of a cycle-registration scheme, a technique described in \cite{RB}. In the context of our problem, this technique can be
described as follows: in addition to the filtrations adjoined to sets $A$ and $B$, we will also observe the filtration adjoined to the set $A\cup B$. The persistence module $PM_k(\mathcal C_{A\cup B})$ can be viewed as a "bigger" module in which persistence modules $PM_k(\mathcal C_A)$ and $PM_k(\mathcal C_B)$ are naturally embedded. More precisely,
these embeddings are morphisms $h^A:PM_k(\mathcal C_A)\Rightarrow PM_k(\mathcal C_{A\cup B})$ and $h^B:PM_k(\mathcal C_B)\Rightarrow PM_k(\mathcal C_{A\cup B})$, such that, at every level $r_i$, mappings $h_{r_i}^A$ and $h_{r_i}^B$ are induced by inclusions. If $\gamma_A$ is a $k-$cycle in the persistence module $PM_k(\mathcal C_A)$ and $\gamma_B$ is a $k-$cycle in $PM_k(\mathcal C_B)$, then these cycles are called $\mathcal C_{A\cup B}-$equivalent cycles (denoted by $\gamma_A\stackrel{\mathcal C_{A\cup B}}{\sim}\gamma_B$), if there are $k-$cycles $\tilde{\gamma}_A\in im(h^A), \tilde{\gamma}_B\in im(h^B)$, such that:
\begin{itemize}
\item Cycles $\gamma_A$ and $\tilde{\gamma}_A$ are born at the same level,
\item Cycles $\gamma_B$ and $\tilde{\gamma}_B$ are born at the same level,
\item Cycles $\tilde{\gamma}_A$ and $\tilde{\gamma}_B$ die at the same level.
\end{itemize}

The notion of $\mathcal C_{A\cup B}-$equivalent cycles is particularly significant in the case when filtrations of the complexes $\mathcal C_A,\mathcal C_B$ and
$\mathcal C_{A\cup B}$ are Morse filtrations. In this case, the first two conditions imply that the cycles $\tilde{\gamma}_A, \tilde{\gamma}_B$ are structurally related to the
cycles $\gamma_A,\gamma_B$, since they appear at the same filtration level. The third condition implies that cycles $\tilde{\gamma}_A$ and $\tilde{\gamma}_B$ are "killed off"
at the same filtration level, which leads to the conclusion that they represent similar homological feature. Consequently, the same conclusion applies to their counterparts, cycles $\gamma_A$ and $\gamma_B$. Also, the following lemma is easily verified.

\begin{lm} Let $\gamma_A$ be a $k-$cycle in the persistence module $PM_k(\mathcal C_A)$ and let $\beta_B,\gamma_B$ be $k-$cycles in the persistence module $PM_k(\mathcal C_B)$
such that $\gamma_A\stackrel{\mathcal C_{A\cup B}}{\sim}\beta_B$ and $\gamma_A\stackrel{\mathcal C_{A\cup B}}{\sim}\gamma_B$. If filtrations of the complexes
$\mathcal C_A,\mathcal C_B$ and $\mathcal C_{A\cup B}$ are Morse filtrations, then $\beta_B=\gamma_B$.
\end{lm}

Therefore, if the filtrations of the complexes $\mathcal C_A,\mathcal C_B$ and $\mathcal C_{A\cup B}$ are Morse filtrations for dimension $k\geqslant 1$, then comparison of the
barcodes $BC_k(A)$ and $BC_k(B)$ can be performed by using matching, which would favor all bars that correspond to the $\mathcal C_{A\cup B}-$equivalent cycles. For those bars in barcodes $BC_k(A)$ and $BC_k(B)$ that cannot be matched in this way, we use the "ordinary" bottleneck distance matching. More details about this "hybrid" matching will be provided in section \ref{subsec:sim.measure}.\\

Unfortunately, the described strategy is troublesome in the case when at least one of the observed filtrations is not a Morse filtration. Taking into account
the discrete nature of the Hamming distance, the possibility that two or more cycles appear or disappear at the same filtration level becomes more and more certain
as the number of strings in the string set increases. In section \ref{sub:string separation} we develop a new technique that will deal with this problem in a
satisfactory way.

\subsection{Separation of simplex radii} \label{sub:string separation}

When analyzing the filtration adjoined to the set $A\subseteq S(n,l)$, its simplices are usually divided into positive (those that mark the birth of a new homology class) and negative (marking the death of such a class). Since $2^{|A|}-1$ simplices must be distributed within at most $l+1$ filtration levels, the scenario in which two or more positive
(or negative) simplices have the same radius is likely to happen. Therefore, there is no guarantee that the filtration adjoined to the set $A$ is a Morse filtration. However, as
we will show, it is possible to construct a set $A'$ of generalized strings such that the filtration adjoined to this set is a Morse filtration. More importantly, this construction
causes strictly controlled "shifts" of bars in the barcode $BC_k(A)$. We begin by giving some definitions.

\begin{de} Closed ball around $x\in S'(n,l)$ with radius $r\geqslant 0$ is the set $B(x,r):=\{y\in S'(n,l):d_{GH}(x,y)\leq r\}$.

$MB(\sigma)=\{B(c,r(\sigma)):c\mbox{ is a center of }\sigma\}$ is the set of miniballs "circumscribed" around the simplex $\sigma$.
\end{de}

Note that a simplex may have more than one center, so that is why we consider the set of miniballs. Notions of this kind were examined in detail in \cite{Z} in the context of the Euclidean space $\mathbb{R}^d$. As usual, the interior of any closed ball $B=B(x,r)$ in the metric space $(S'(n,l),d_{GH})$ is ${\rm int}B:=\{y\in S'(n,l):d_{GH}(x,y)<r\}$ and its boundary is ${\rm bd}B:=\{y\in S'(n,l):d_{GH}(x,y)=r\}$.

\begin{de}
For a finite subset $A\subseteq S'(n,l)$, $G\subseteq A$ is called a set of generators if there is $B\in MB(A)$ such that $G\subseteq{\rm bd}B$ and $A\setminus G\subseteq{\rm int}B$.
\end{de}

\begin{lm}\label{minimal}
Every finite $A\subseteq S'(n,l)$ has a minimal set of generators.
\end{lm}

\dokaz Let $G_1$ and $G_2$ be two sets of generators for some finite $A\subseteq S'(n,l)$. Let $B_1$ and $B_2$ be the corresponding miniballs with radius $r$, $c_1$ and $c_2$ their centers, and let $G=G_1\cap G_2$. If $c[i](j):=\frac{c_1[i](j)+c_2[i](j)}2$, then $A\subseteq B(c,r)$. Indeed, for every $x\in A$,

\begin{eqnarray} \label{formula4}
d_{GH}(c,x) &=& \sum_{i=1}^l\left(1-\sum_{j=1}^n\min\left\{\frac{c_1[i](j)+c_2[i](j)}2,x[i](j)\right\}\right) \nonumber\\
 &\leq & \sum_{i=1}^l\left(1-\sum_{j=1}^n\frac 12\big(\min\{c_1[i](j),x[i](j)\}+\min\{c_2[i](j),x[i](j)\}\big)\right)\nonumber\\
 &=& \frac 12\sum_{i=1}^l\big(1-\sum_{j=1}^n\min\{c_1[i](j),x[i](j)\}\big)\nonumber\\
 &+&\frac{1}{2}\sum_{i=1}^l(1-\sum_{j=1}^n\min\{c_2[i](j),x[i](j)\})\nonumber\\
 &=& \frac 12\big(d_{GH}(c_1,x)+d_{GH}(c_2,x)\big)\leq r.
\end{eqnarray}

Note that the inequality given in (\ref{formula4}) can be equality only for $x\in G$. Hence
the assumption that $G_1$ and $G_2$ are disjoint sets would lead to the conclusion that there is $r'<r$ such that $d_{GH}(c,x)\leqslant r'$ holds for every
$x\in A$, which is impossible because $r(A)=r$. So $G$ is nonempty and "generates" another miniball of radius $r$ circumscribed around $A$. Hence, the intersection of sets of generators contains another set of generators, which means
that the intersection of them all is the minimal set of generators.\kraj

\begin{ex}
Let $s_1=1111$, $s_2=2222$, $t=1222$, $u=1212$, $c_1=1122$, $c_2=2211$ and let $c_3=c_3[1]c_3[2]c_3[3]c_3[4]\in S'(2,4)$ be given by $c_3[1]=c_3[2]=\left(\begin{array}{cc}
1 & 2\\
\frac 12 & \frac 12
\end{array}\right)$, $c_3[3]=\left(\begin{array}{cc}
1 & 2\\
1 & 0
\end{array}\right)$ and $c_3[4]=\left(\begin{array}{cc}
1 & 2\\
0 & 1
\end{array}\right)$. Then each of $c_1,c_2,c_3$ is a center for the simplex $\sigma=\{s_1,s_2\}$ and the radii of the corresponding miniballs are 2. However, since $t\in B(c_1,2)\setminus B(c_2,2)$, the first of these two miniballs is also circumscribed around $\tau=\sigma\cup\{t\}$, but the second is not. Thus, $\sigma$ is the minimal set of generators for both $\sigma$ and $\tau$. For $\theta=\sigma\cup\{u\}$, $\theta$ itself is a set of generators (since all vertices of $\theta$ lie on the boundary of $B(c_1,2)$), but $\sigma$ is the minimal one: $s_1,s_2\in{\rm bd}B(u,2)$ while $u\in{\rm int}B(u,2)$.
\end{ex}

\begin{de}
Define a binary relation $\approx$ on finite subsets of $S'(n,l)$ as follows: $A\approx B$ if $A$ and $B$ have the same minimal set of generators.
\end{de}

Clearly, $\approx$ is an equivalence relation. It will turn out that, for a given filtration, the simplices that can not be separated (at least not by the method described below) are exactly those that are in the same $\approx$-equivalence class.

For $s\in S'(n,l)$ and $k<l$, $s\upharpoonright\bN_k\in S'(n,k)$ denotes the generalized string consisting of the first $k$ elements of $s$.  For given $\sigma\in \mathcal C_{A}$, let $C(\sigma)$ be the set of centers of miniballs circumscribed around the minimal set of generators $G$ and let $D(\sigma,u)=\min\{d_{GH}(c,u):c\in C(\sigma)\}$.

\begin{lm}\label{pomseparation}
Let $A\subseteq S'(n,l)$ be a finite set of generalized strings and let $\sigma_1,\sigma_2\in \mathcal C_{A}$ be simplices such that $r(\sigma_1)=r(\sigma_2)=r_0$ and $\sigma_1\not\approx\sigma_2$. Also, let $j\in\mathbb{N}$ be arbitrary. Then, there are set $B\subseteq S'(n,l+1)$, a vertex $z\in A$ and a bijection $f:A\rightarrow B$ such that $r(f[\sigma_1])=r(\sigma_1)$, $r(f[\sigma_2])>r(\sigma_2)$ and
\begin{equation}\label{eqsep}
r(\tau)\leq r(f[\tau])\leq r(\tau)+\frac 1j
\end{equation}
for all $\tau\in \mathcal C_{A}$. Furthermore:

(i) if $\frac 1j<\min\big(\{|r(\tau)-r(\sigma)|:\sigma,\tau\in \mathcal C_{A}\}\setminus\{0\}\big)$, then $r(\sigma)<r(\tau)$ implies $r(f[\sigma])<r(f[\tau])$, for all $\sigma,\tau\in \mathcal C_{A}$;

(ii) if $\sigma\in \mathcal C_{A}$, $G$ is the minimal set of generators for $\sigma$, $z\notin G$ and $\frac 1j<r(\sigma)-D(\sigma,z)$, then $f[G]$ is the minimal set of generators for $f[\sigma]$.
\end{lm}

\dokaz Let $G_1$ and $G_2$ be the minimal sets of generators for $\sigma_1$ and $\sigma_2$. The condition $\sigma_1\not\approx\sigma_2$ means that, say, $G_2\not\subseteq G_1$.
So, we can pick a generalized string $z\in G_2\setminus G_1$. Now define, for any $s=a_1a_2\dots a_l\in A$, $f(s)$ as follows: $f(s)=a_1a_2\dots a_la_{l+1}$, where:

- for $s\neq z$, let $a_{l+1}(1):=1$ and $a_{l+1}(i):=0$, for $i>1$, and

- for $s=z$, let $a_{l+1}(1):=1-\frac 1j$, $a_{l+1}(2):=\frac 1j$ and $a_{l+1}(i):=0$, for $i>2$.

Now, if $c=b_1b_2\dots b_l$ is the center of the miniball $B(c,r_0)$ circumscribed around $\sigma_1$, then $c':=b_1b_2\dots b_lb_{l+1}$ (where $b_{l+1}(1):=1$ and $b_{l+1}(i)=0$ for $i>1$) is the center of the closed ball with radius $r_0$ containing $f[\sigma_1]$, and so $r(f[\sigma_1])=r_0$.

In a similar way, we see that (\ref{eqsep}) holds for any $\tau\in \mathcal C_{A}$.

For any two generalized strings $c\in S'(n,l+1)$ and $y\in\sigma_2$, we have
\begin{eqnarray*}
d_{GH}(c,f(y)) &=& \sum_{i=1}^{l+1}\left(1-\sum_{j=1}^n\min\left\{c[i](j),f(y)[i](j)\right\}\right)\\
 &=& \sum_{i=1}^l\left(1-\sum_{j=1}^n\min\{c[i](j),f(y)[i](j)\}\right)\\
 &+&\big(1-\sum_{j=1}^n\min\{c[l+1](j),f(y)[l+1](j)\}\big)\\
 &=& d_{GH}(c\upharpoonright\mathbb{N}_l,y)+\left(1-\sum_{j=1}^n\min\left\{c[l+1](j),f(y)[l+1](j)\right\}\right).
\end{eqnarray*}
If we assume that, for some $c_0$, $d_{GH}(c_0,f(y))\leq r_0$ for every $y\in\sigma_2$, it follows that $d_{GH}(c_0\upharpoonright\mathbb{N}_l,y)\leq r_0$ for every $y\in G_2$, so $c_0\upharpoonright\mathbb{N}_l$ must be a center of a miniball of $\sigma_2$. However, for each such $c_0$ we have $c_0[l+1]\neq z[l+1]$, so  $1-\sum_{j=1}^n\min\{c_0[l+1](j),f(z)[l+1](j)\}>0$ and consequently $d_{GH}(c_0,f(z))>r_0$. Hence, $r(f[\sigma_2])$ must be greater than $r_0$.

(i) follows easily from (\ref{eqsep}). Finally, for (ii), the condition $\frac 1j<r(\sigma)-D(\sigma,z)$ guarantees that, since $z$ was in the interior of some miniball $B(c,r)$ circumscribed around $G$, then $f(z)$ belongs to interior of at least one miniball (namely $B(f(c),r)$) circumscribed around $f[G]$.\kraj

It should be noted that the bijection described in the previous lemma rightshifts levels of the persistence module $PM_k\left(\mathcal C_{A}\right)$ for at most $\frac{1}{j}$.
This fact, together with The Stability Theorem, implies that the bottleneck distance between barcodes $BC_k(A)$ and $BC_k(B)$ is less than or equal to
$\frac{1}{j}$.\\

After one application of the previous lemma, it is still possible that there are non $\approx$-equivalent simplices with the same radius in the full complex $\mathcal C_{B}$. In order to "separate" radii of those simplices, we will successively continue to apply this lemma, with the appropriate choice of $\frac{1}{j}$, which will ensure that, in each of these steps,
the radii of the simplices that we separated earlier do not become equal again.

\begin{te}\label{separacija}
Let $A\subseteq S'(n,l)$ be such that $|A|=m$ and let $\varepsilon>0$ be given. Then there are $Sep(A)\subseteq S'(n,l+m')$, for some $m'\leq m$, and a bijection
$g:A\rightarrow Sep(A)$ such that:

(i) $r(g[\sigma])\neq r(g[\tau])$ for all $\sigma,\tau\in \mathcal C_{A}$ such that $\sigma\not\approx\tau$, and

(ii) $0\leq r(g[\sigma])-r(\sigma)<\varepsilon$ for all $\sigma\in \mathcal C_{A}$.
\end{te}

\dokaz We use Lemma \ref{pomseparation} several times, each time separating two simplices and changing the radii of others for sufficiently small amounts. First, let $A_0:=A$
and let $\sigma_1,\sigma_2\in \mathcal C_{A_0}$ be simplices such that $\sigma_1\not\approx\sigma_2$, $r(\sigma_1)=r(\sigma_2)$. 
Choose $z$ from the minimal set of generators of, say, $\sigma_1$ as in Lemma \ref{pomseparation}, and let $j_1\in\mathbb{N}$ be such that 
\begin{equation*}
\resizebox{\hsize}{!}{$\frac 1{j_1}<\min\left\{\frac\varepsilon 2,\min(\{r(\sigma)-D(\sigma,z):\sigma\in \mathcal C_{A}\}\cap\mathbb{R}^+),\min(\{|r(\tau)-r(\sigma)|:\sigma,\tau\in \mathcal C_{A_0}\}\setminus\{0\})\right\}$}.
\end{equation*}
We obtain $A_1\subseteq S'(n,l+1)$ and a bijection $f_1:A_0\rightarrow A_1$, such that $r(f_1[\sigma_1])<r(f_1[\sigma_2])$ and $r(f_1[\sigma])<r(f_1[\tau])$, whenever $r(\sigma)<r(\tau)$, for $\sigma,\tau\in \mathcal C_{A_0}$. Now, we repeat the process, using some $j_i$ satisfying  $\frac 1{j_i}<\frac{\varepsilon}{2^i}$, obtaining sets
$A_2,A_3$, $\dots,A_{m'}$, so that in $Sep(A):=A_{m'}$ all simplices that are not $\approx$-equivalent have different radii.
This proves (i). Note that the condition (ii) of Lemma \ref{pomseparation} implies that, if $\sigma\not\approx\tau$, then $f[\sigma]\not\approx f[\tau]$.

In the end, we take $g:=f_{m'}\circ\dots\circ f_2\circ f_1$. Clearly, $0\leq r(g[\sigma])-r(\sigma)\leq\frac 1{j_1}+\frac 1{j_2}+\dots+\frac 1{j_{m'}}<\varepsilon$, for every $\sigma\in \mathcal C_{A}$, which proves (ii). Also, $m'$ will be no larger than $m$ since each vertex $z$ needs to be "moved" at most once (after the moving it can not be a member of another difference $G_2\setminus G_1$ of sets of generators of simplices with the same radius).\kraj

In particular, the condition (ii) in the previous theorem shows that "new" bars (appearing in the barcode $BC_k(Sep(A))$, but not in the barcode of $BC_k(A)$) are of length
less than $\varepsilon$, and the length of each "old" bar of the barcode $BC_k(A)$ has changed for less than $\varepsilon$. Also, we can see that only equivalent simplices can eventually have a same radius in the full complex $\mathcal C_{Sep(A)}$. So  let us show that such equivalence classes of simplices do not affect the barcode $BC_k(Sep(A))$.

\begin{te} \label{ekv.simpl.}
Let $E$ be a  $\,\approx-$equivalence class with at least two elements, and let $r_0>0$ be the radius of all $\sigma\in E$. Then the appearance of simplices from $E$ does not affect the barcode; more precisely: persistence modules $PM_k(\mathcal C_{Sep(A)}^{(r_0)}\setminus E)$ and
$PM_k(\mathcal C_{Sep(A)}^{(r_0)})$ are equal for each dimension $k>1$.
\end{te}

\dokaz Since $\mathcal C_{Sep(A)}$ is obtained as a result of applying Theorem \ref{separacija}, the only simplices with radius $r_0$ in this complex are those in $E$. Let $G$ be the common minimal set of generators for $\sigma\in E$. This means that $E$ consists of all simplices $\sigma$ such that $G\subseteq\sigma$ and $\sigma\setminus G\subseteq{\rm int}B(c,r_0)$, where $B(c,r_0)$ is the ball circumscribed around $G$. Let $x\in {\rm int}B(c,r_0)\setminus G$ be an arbitrary vertex belonging to some of these simplices. All simplices in $E$ can be divided into pairs $(\sigma,\sigma\cup\{x\})$, where $x\notin\sigma$. Let $\langle(\sigma_i,\tau_i):i<d\rangle$ be an enumeration of all such pairs, such that $|\sigma_i|\leq|\sigma_j|$ for $i<j$. Now fix a small enough $\delta$, and let us examine the effect of $E$ on the bar code by "pretending" that the simplices from $E$ appear one by one in order of indices $i$, for example that $r(\sigma_i)=r_0+2i\delta$ and $r(\tau_i)=r_0+(2i+1)\delta$. For this new filtration (call it ${\mathcal K}$) we have ${\mathcal K}^{(r_0)}=\mathcal C_{Sep(A)}^{(r_0)}\setminus E$ and ${\mathcal K}^{(r_0+(2d+1)\delta)}=\mathcal C_{Sep(A)}^{(r_0)}$.

Now fix some $i$ and let $m:=|\sigma_i|$. All $m$-element subsets of $\tau_i=\sigma_i\cup\{x\}$ except $\sigma_i$ have radii less than $r_0+2i\delta$. Indeed, any such subset either does not contain $G$ (in which case their radius is smaller than $r_0$: if $G'$ is a minimal set of generators of such a $\sigma$, then by the proof of Lemma \ref{minimal} $G\cap G'$ also contains a set of generators, so $G'\subset G$), or is of the form $\tau_j$, for some $j<i$. Hence, $\sigma_i$ is a positive simplex, marking the birth of an $m$-dimensional homology class, and $\tau_i$ is the negative simplex killing that same class. Thus, returning to the situation in which all the simplices in $E$ appear simultaneously, their overall effect on the
barcode is none.
\kraj

\subsection{A new string similarity measure} \label{subsec:sim.measure}

We have made all the necessary preparations to introduce a new measure of similarity between two sets of strings. \\

Let $A,B\subseteq S(n,l)$ be two sets of strings, such that $|A|=|B|=m\geqslant 2$. For each dimension $k\geqslant 0$, we will propose a new hybrid matching of $k$-dimensional
bars and define the distance $d_k$ between appropriate barcodes. In this hybrid matching, the priority will be to match bars that correspond to equivalent cycles.\\

For $k=0$, we have already established that both of the barcodes $BC_0(A)$ nad $BC_0(B)$ contain $m$ bars, where $m-1$ are finite-length bars (all having
$0$ as a left endpoint) and one bar is the infinite-length bar. If $0<l_1^A\leqslant l_2^A\leqslant\dots\leqslant l_{m-1}^A$  and
$0<l_1^B\leqslant l_2^B\leqslant\dots\leqslant l_{m-1}^B$ are lengths of finite-length bars, then we can match bars $[0,l_i^A)$ and $[0,l_i^B)$ and define the distance $\displaystyle d_0(A,B):=\max_{i\in\{1,2,\dots,m-1\}}|l_i^A-l_i^B|$.\\

For a dimension $k\geqslant 1$, we use our simplices radii separation technique to get $m-$element sets $Sep(A)$ and $Sep(B)$ of generalized strings. If there are no bars in either of the barcodes $BC_k(Sep(A))$ and $BC_k(Sep(B))$, we set $d_k(A,B):=0$. Otherwise, we apply simplices radii separation technique one more time to get the set $Sep(A\cup B)$
of generalized strings.  Note that the separation in $A\cup B$ can be performed by including the steps of the separation in both of $A$ and $B$, so that
$Sep(A)\subseteq Sep(A\cup B)$ and $Sep(B)\subseteq Sep(A\cup B)$.
In this way, we ensure that $\mathcal C_{Sep(A)}$, $\mathcal C_{Sep(B)}$ and $\mathcal C_{Sep(A\cup B)}$ are Morse filtrations. Next, we look for a potential $\mathcal C_{Sep(A\cup B)}-$equivalent cycles and match their corresponding bars. For  bars in barcodes $BC_k(Sep(A))$ and $BC_k(Sep(B))$ which are not matched in this way, we use bottleneck distance matching. More precisely, if $BC_k'(Sep(A))\subseteq BC_k(Sep(A))$ and $BC_k'(Sep(B))\subseteq BC_k(Sep(B))$ denote collections of all bars without any
$\mathcal C_{Sep(A\cup B)}-$equivalent counterpart, then we can define
\begin{equation} \label{d_k}
d_k(A,B):=\sum_{\gamma_{1}{\sim}\gamma_{2}} d_{BOT}\big(\{l(\gamma_{1})\},\{l(\gamma_{2})\}\big)+d_{BOT}\big(BC_k'(Sep(A)),BC_k'(Sep(B))\big),
\end{equation}
where the first sum is taken over all pairs of $\mathcal C_{Sep(A\cup B)}-$equivalent $k-$cycles $\gamma_{1},\gamma_{2}$, and
$l(\gamma_{1})\in BC_k(Sep(A)),l(\gamma_{2})\in BC_k(Sep(B))$ are bars corresponding to these cycles. Of course, in the case when there are no $\mathcal C_{Sep(A\cup B)}-$equivalent cycles, we have $d_k(A,B)=d_{BOT}\big(BC_k(Sep(A)),BC_k(Sep(B))\big)$. The comparison of barcode lines of $BC_k(Sep(A))$ and $BC_k(Sep(B))$ is justified
by a fact that, for every $\varepsilon>0$,  sets $Sep(A)$ and $Sep(B)$ can be chosen such that
\begin{align*}
&d_{BOT}\big(BC_k(A),BC_k(B)\big)\leqslant \underbrace{d_{BOT}\big(BC_k(A),BC_k(Sep(A))\big)}_{\leqslant \frac{\varepsilon}{2}}\\
&+d_{BOT}\big(BC_k(Sep(A)),BC_k(Sep(B))\big)+\underbrace{d_{BOT}\big(BC_k(Sep(B),BC_k(B)\big)}_{\leqslant \frac{\varepsilon}{2}}\\
&\leqslant d_{BOT}\big(BC_k(Sep(A)),BC_k(Sep(B))\big)+\varepsilon.
\end{align*}
Let $k_0\geqslant 0$ be a minimal dimension with property that $BC_k(Sep(A))=\emptyset=BC_k(Sep(B))$, for every $k>k_0$. We define a new distance measure between sets
$A,B\subseteq S(n,l)$ of the same cardinality by
$$ d_{new}(A,B):=\sum_{k=0}^{k_0} \frac{2^k}{2^{k_0+1}-1}\cdot d_k(A,B).  $$
Weights $\frac{2^k}{2^{k_0+1}-1}$, $0\leqslant k\leqslant k_0$, are assigned in order to prioritize differences in a homology features of sets $A$ and $B$, in the favor of those discrepancies that are manifested
in higher dimensions. The distance $d_{new}$ has the stability property, since every distance $d_k$ is defined via the bottleneck distance between appropriate sets of barcodes.

\section{Conclusions and future work}

In many disciplines, including information theory, coding theory, cryptography, and bioinformatics, strings are used to encode finite sequential data types. Examination of
measures of similarity between two sets of strings is an ongoing investigation of various patterns that would enable the comparison of these sets. In this paper, we use
the tools from persistence homology in order to quantify the similarity of "connectivity issues" of various dimensions that may exist for given sets of strings. This is
accomplished by constructing the new measure $d_{new}$ based on the newly proposed hybrid matching, whose main property is giving priority to matching barcode lines of the corresponding equivalent cycles. The applicability of our hybrid matching is heavily dependent on an assumption that all involved filtrations are Morse filtrations. To
fullfill this requirement, we develope the separation of simplex radii technique, which we introduce in Lemma \ref{pomseparation} and Theorem \ref{separacija}.
Also, we identify a notion of $\approx-$equivalent simplices (simplices with the same minimal set of generators) and become aware of their interesting property that they are not affecting the structure of barcode lines. This is stated in Theorem \ref{ekv.simpl.}. It is important to point out that this result can be viewed in a broader context that does not necessarily include the analysis of string similarity measures.\\

More work will be needed to construct efficient algorithms for conducting the ideas of this paper. In particular: (1) calculating the radius and set of centers of a given finite subset of $S'(n,l)$, (2) choosing pairs $(\sigma_1,\sigma_2)$ in Theorem \ref{separacija} to minimize the number of steps, and hence the dimension of the obtained space and (3) if possible, performing the process of the theorem so that we do not need to calculate radii from the beginning each time, but to get them from the previous values of radii.\\

For future work, the authors would like to investigate a potential sufficient condition under which assumption of the existence of a filtration isomorphism
would guarantee existence of a $d_H-$isomorphism between appropriate sets of strings. Also, we would like to use the methodology presented in this paper for the purpose of
developing string similarity measures based on some other string metrics. More concretely, we would like to investigate string similarity measures based on the longest common subsequence ($LCS$) metric.
It would be useful to find an analogy for the separation of simplex radii technique in this case. Also, it would be very nice to appraise the role of $\approx-$equivalent simplices as some sort of "neutral" packs of simplices in the general \v Cech filtration setup.

\end{document}